\documentclass[11pt]{article}
\usepackage{amsmath,latexsym,epsf,mathabx}
\begin{document}

\newcommand{\nc}{\newcommand}
\def\PP#1#2#3{{\mathrm{Pres}}^{#1}_{#2}{#3}\setcounter{equation}{0}}
\def\ns{$n$-star}\setcounter{equation}{0}
\def\nt{$n$-tilting}\setcounter{equation}{0}
\def\Ht#1#2#3{{{\mathrm{Hom}}_{#1}({#2},{#3})}\setcounter{equation}{0}}
\def\qp#1{{${(#1)}$-quasi-projective}\setcounter{equation}{0}}
\def\mr#1{{{\mathrm{#1}}}\setcounter{equation}{0}}
\def\mc#1{{{\mathcal{#1}}}\setcounter{equation}{0}}
\def\HD{\mr{Hom}_{\mc{D}}}
\def\HC{\mr{Hom}_{\mc{C}}}
\def\AdT{\mr{Add}_{\mc{T}}}
\def\adT{\mr{add}_{\mc{T}}}
\def\Kb{\mc{K}^b(\mr{Proj}R)}
\def\kb{\mc{K}^b(\mc{P}_R)}
\def\AdpC{\mr{Adp}_{\mc{C}}}
\newtheorem{Th}{Theorem}[section]
\newtheorem{Def}[Th]{Definition}
\newtheorem{Lem}[Th]{Lemma}
\newtheorem{Pro}[Th]{Proposition}
\newtheorem{Cor}[Th]{Corollary}
\newtheorem{Rem}[Th]{Remark}
\newtheorem{Exm}[Th]{Example}
\newtheorem{Sc}[Th]{}
\def\Pf#1{{\noindent\bf Proof}.\setcounter{equation}{0}}
\def\>#1{{ $\Rightarrow$ }\setcounter{equation}{0}}
\def\<>#1{{ $\Leftrightarrow$ }\setcounter{equation}{0}}
\def\bskip#1{{ \vskip 20pt }\setcounter{equation}{0}}
\def\sskip#1{{ \vskip 5pt }\setcounter{equation}{0}}
\def\mskip#1{{ \vskip 10pt }\setcounter{equation}{0}}
\def\bg#1{\begin{#1}\setcounter{equation}{0}}
\def\ed#1{\end{#1}\setcounter{equation}{0}}
\def\KET{T^{^F\bot}\setcounter{equation}{0}}
\def\KEC{C^{\bot}\setcounter{equation}{0}}


\title{\bf  Quasi-cotilting modules and torsion-free classes\thanks{Supported by the National Science Foundation of China (Grant Nos. 11371196 and 11171149) and the National Science Foundation for Distinguished Young Scholars of Jiangsu Province (Grant No. BK2012044).}}

\smallskip
\author{\small Peiyu Zhang and Jiaqun Wei\thanks{Correspondent auther}
\\
\small E-mail:~zhangpeiyu2011@163.com,~weijiaqun@njnu.enu.cn\\
\small Institute of Mathematics, School of Mathematics Sciences\\
\small Nanjing Normal University, Nanjing \rm210023 China}
\date{}
\maketitle
\baselineskip 15pt
%
%
\begin{abstract}
\vskip 10pt%
We prove that all quasi-cotilting modules are pure-injective and cofinendo.  It follows that the class $\mathrm{Cogen}M$ is always a covering class whenever $M$ is a quasi-cotilting module. Some characterizations of quasi-cotilting modules are given. As a main result, we prove that there is a bijective correspondence between the equivalent classes of quasi-cotilting  modules and torsion-free covering classes.

\mskip\

\noindent MSC 2010: 16D90 18E40 16E30 16S90


\noindent {\it Keywords}: quasi-cotilting module, torion-free class, covering class, pure-injective module, cofinendo module

\end{abstract}
%
\vskip 30pt

\section{Introduction and Preliminaries}

%
%
%
\hskip 18pt

Quasi-tilting modules were introduced by Colpi, D’este, and Tonolo in \cite{RCDT} in the study of contexts of $\ast$-modules and tilting modules, where it was shown that these modules are closely relative to torsion theory counter equivalences. Recently, quasi-tilting modules turn out to be new interesting after the study of support $\tau$-tilting modules by Adachi, Iyama and Reiten \cite{AIR}. The second author proved that a finitely generated module over an artin algebra is support $\tau$-tilting if and only if it is quasi-tilting \cite{W14}. Moreover, Angeleri-H\"{u}gel, Marks and Vit\'{o}ria \cite{AMV} proved that finendo quasi-tilting modules are also closely related to silting modules (which is a generalizations of support $\tau$-tilting modules in general rings).

It is natural to consider the dual of quasi-tilting modules, i.e., quasi-cotilting modules. As cotilting modules possess their own interesting properties not dual to tilting modules, we show in this paper that quasi-cotilting modules are not only to be the dual of quasi-tilting modules and they have their own interesting properties too. We prove that all quasi-cotilting modules are pure-injective and cofinendo. These are clearly new important properties of quasi-cotilting modules. Note that not all quasi-tilting modules are finendo. As a corollary, we obtain that the class $\mathrm{Cogen}M$ is a covering class whenever $M$ is quasi-cotilting. We give variant characterizations of quasi-cotilting modules and cotilting modules. As the main result, we prove that there is a bijection between the equivalent classes of quasi-cotilting  modules and torsion-free covering classes.

\vskip 10pt

Throughout this paper, $R$ is always an associative ring with identity and subcategories are always full and closed under isomorphisms.
We denote by $R$-Mod the category of all left $R$-modules. By $\mathrm{Proj}R$ we denote
the class  of all projective $R$-module.
%


Let $M\in R$-Mod, we use following notations throughout this paper.

\bg{verse}
$\mathrm{Adp} M:=\{N\in R$-Mod $\mid$ there is a module $L$ such that $N\oplus L=M^{X}$ for some $X\}$;

$\mathrm{Cogen} M:=\{N\in R$-Mod $\mid$ there is an exact sequence $0\to N\to M_{0}$ with $M_{0} \in \mathrm{Adp} M\}$;

$\mathrm{Copres} M:=\{N\in R$-Mod $\mid$ there is an exact sequence $0\to N\to M_{0}\to M_{1}$ with $M_{0}$,
$M_{1} \in \mathrm{Adp} M\}$;

$^{\bot_{1}}M:=\{N\in R$-Mod $\mid$ $\mathrm{Ext}^{1}_R(N,M)=0 \}$;

$M^{\bot_{1}}:=\{N\in R$-Mod $\mid$ $\mathrm{Ext}^{1}_R(M,N)=0 \}$;

$^{\circ}M:=\{N\in R$-Mod $\mid \mathrm{Hom}_R(N,M)=0 \}$;

$M^{\circ}:=\{N\in R$-Mod $\mid \mathrm{Hom}_R(M,N)=0 \}$.
\ed{verse}

Note that $\mathrm{Cogen}M$ is clearly closed under submodules, direct products and direct sums.

Let $\mathcal{T}$ be a class of $R$-modules, we denoted by $\mathrm{Fac}(\mathcal{T})$ the classes formes by the factor modules of all
modules in $\mathcal{T}$.
An $R$-module $M \in \mathcal{T}$ is called Ext-injective in $\mathcal{T}$
if $\mathcal{T} \subseteq$ $^{\bot_{1}}M$.





Let $\mathcal{X}$ and $\mathcal{Y}$ be two subcategories of $R$-Mod. The pair ($\mathcal{X}$,$\mathcal{Y}$) is
said to be a torsion pair if it satisfies the following three condition (1) $\mathrm{Hom} (\mathcal{X},\mathcal{Y})=0$; (2) if $\mathrm{Hom}(M,\mathcal{Y})=0$,
then $M \in \mathcal{X}$; (3) if $\mathrm{Hom}(\mathcal{X},N)=0$, then $N \in \mathcal{Y}$.
In the case, $\mathcal{X}$ is called a torsion class and $\mathcal{Y}$ is called a torsion-free class.
Note that a subcategory $\mathcal{X}$ of $R$-Mod is torsion-free if and only if $\mathcal{X}$ is closed under direct products,
submodules and extensions, see \cite{SED}.
\bg{Lem}\label{IARM}%
If an R-module M is Ext-injective in $\mathrm{Cogen}M$,
then $(^{\circ}M, \mathrm{Cogen}M)$ is a torsion pair.

\ed{Lem}
\Pf. It is easy to verify that  $^{\circ}M=$$^{\circ}(\mathrm{Cogen}M)$.
We only need to prove that  $(^{\circ}M)^{\circ}=\mathrm{Cogen}M$.
If $T\in \mathrm{Cogen}M$, then there is an injective homomorphism $i$: $T\to M^{X}$
for some set $X$. For any $N\in$ $^{\circ}M$ and $f\in \mathrm{Hom}_R(N,T)$, then $\mathrm{Hom}_R(N,M^{X})=
(\mathrm{Hom}_R(N,M))^{X}=0$, hence $if=0$. Then  $f=0$ since $i$ is injective.
So $\mathrm{Cogen}M \subseteq (^{\circ}M)^{\circ}$.

For any $T\in(^{\circ}M)^{\circ}$, consider the evaluation map $\alpha$: $T\to M^{\mathrm{Hom}_R(T,M)}$
with $K=\ker \alpha$. Take canonical resolution of $\alpha$, i.e. $\alpha=i \pi$
with $\pi$: $T\to \mathrm{Im} \alpha$ and $i$: $\mathrm{Im}\alpha \to M^{\mathrm{Hom}_R(T,M)}$.
Clearly, $\mathrm{Hom}_R(\pi,M)$ is surjective from the definition of $\alpha$.
Applying the functor $\mathrm{Hom}_R(-,M)$ to the exact sequence $0\to K\to T\to Im \alpha\to0$,
we have an exact sequence
$$0\to \mathrm{Hom}_R(Im \alpha,M)\to \mathrm{Hom}_R(T,M)\to \mathrm{Hom}_R(K,M)\to \mathrm{Ext}^{1}_R(Im \alpha,M).$$
Since $Im \alpha\in \mathrm{Cogen}M\subseteq$$^{\bot_{1}}M$, $\mathrm{Ext}^{1}_R(Im \alpha,M)$$=0$.
As $\mathrm{Hom}_R(\pi,M)$ is surjective by above discussion, $\mathrm{Hom}_R(K,M)$ $=0$, and hence, $K\in$ $^{\circ}M$.
Since $T\in$$(^{\circ}M)^{\circ}$, we have that $\mathrm{Hom}_R(K,T)=0$ and $\alpha$
is injective. Thus $K=0$, i.e. $T\in \mathrm{Cogen}(M)$.
\ \hfill $\Box$


\bg{Def}\label{}$\mathrm{\cite[Definition 1.4]{HTTT}}$%
\ $(1)$ Let $Q$ be an injective cogenerator of $R$-Mod. A module $M$ is called $Q$-cofinendo
if there exist a cardinal $\gamma$ and a map $f$: $M^{\gamma}\to Q $ such that for any
cardinal $\alpha$, all maps $M^{\alpha}\to Q$ factor through $f$.

$(2)$ A module $M$ is cofinendo if there is an injective cogenerator $Q$ of $R$-Mod
such that $M$ is
$Q$-cofinendo.
\ed{Def}

Let $\mathcal{T}$ be a class of $R$-modules and $M$ be an $R$-module. Then $f$: $X\to M$ with $X\in \mathcal{T}$
is a $\mathcal{T}$-precover of
$M$ provided that $\mathrm{Hom}(Y,f)$ is surjective for any $Y\in\mathcal{T}$.
A $\mathcal{T}$-precover of $M$ is called $\mathcal{T}$-cover of $M$ if any $g$: $X\to X$ such that $f=fg$ must be an isomorphism.
A class $\mathcal{T}$ of $R$-modules is said to be a precover class (cover class) provided that each module
has a $\mathcal{T}$-precover ($\mathcal{T}$-cover).

\bg{Lem}\label{TFAEF}$\mathrm{\cite[Proposition 1.6]{HTTT}}$%
\ The following are equivalent for a module $M$:

$(1)$ $M$ is cofinendo;

$(2)$ there is an $\mathrm{Adp} M$-precover of an injective cogenerator Q of $R$-Mod;

$(3)$ $\mathrm{Cogen} M$ is a precover class.

\ed{Lem}

\mskip\


\section{Quasi-cotilting modules}
\mskip\

\hskip 25pt In this section, we introduce notion of quasi-cotilting modules and
give some characterization of  quasi-cotilting modules. In particular, we prove
that all quasi-cotilting modules are pure injective and cofinendo.

\bg{Def}\label{}%
$(1)$ An $R$-module $M$ is said to be a costar module if $\mathrm{Cogen}M=\mathrm{Copres}M$
and $\mathrm{Hom}_R(-,M)$ preserves exactness of  any short exact sequence in $\mathrm{Cogen}M$.

$(2)$ An $R$-module $M$ is called a quasi-cotilting module, if it is  a costar module and
M is Ext-injective in $\mathrm{Cogen} M$.
\ed{Def}

{\noindent\bf Remark}  (1) Costar modules defined above had been studied in \cite{LHX}. Moreover, their general version, i.e., $n$-costar modules, were studied by He \cite{Hd} and Yao and Chen \cite{YC} respectively.

(2) Colby and Fuller \cite{CF} had defined another notion of costar modules which can be viewed as a special case of the above-defined costar modules.

\vskip 10pt
For convenience, we say that a short exact sequence $0\to A\to B\to C\to 0$ is $\mathrm{Hom}_R(-,M)$-exact (($M\bigotimes_{R}-$, resp.)-exact) if
the functor $\mathrm{Hom}_R(-,M)$ ($M\bigotimes_{R}-$, resp.) preserves exactness of this exact sequence.

The following result is well-known.

\bg{Lem}\label{direct}%
Suppose that two short exact sequences $0\to A\to B\to C\to 0$ and $0\to A\to B^{'}\to C^{'}\to 0$
are $\mathrm{Hom}(-,M)$-exact with $B,B'\in \mathrm{Apd}M$. Then $B\bigoplus C^{'}\cong B^{'}\bigoplus C$

\ed{Lem}

\Pf. It is dual to Lemma 2.2 in \cite{WJQ}.
\ \hfill $\Box$

\mskip\

The following result presents a useful
property of costar modules.

\bg{Lem}\label{LMBAC}%
Let $M$ be a co-star module. Suppose that the short exact sequence $0\to A\to B\to C\to 0$ is
$\mathrm{Hom}_R(-,M)$-exact and $A \in \mathrm{Cogen}M$, then $B \in \mathrm{Cogen}M$ if and
only if $C \in \mathrm{Cogen}M$.

\ed{Lem}

\Pf. $\Leftarrow$ Since $A$ and $C$ are in $\mathrm{Cogen}M$, we have two monomorphisms $f$: $A\to M_{A}$
and $g$: $C\to M_{C}$ with $M_{A}$ and $M_{C}$ in $\mathrm{Adp} M$. We consider the
following diagram:

 \setlength{\unitlength}{0.09in}
 \begin{picture}(50,10)

 \put(14,1){\makebox(0,0)[c]{$0$}}
                              \put(16,1){\vector(1,0){3}}
 \put(21,1){\makebox(0,0)[c]{$M_{A}$}}
                             \put(23,1){\vector(1,0){3}}
 \put(31,1){\makebox(0,0)[c]{$M_{A}\bigoplus M_{C}$}}
                             \put(36,1){\vector(1,0){3}}
 \put(41,1){\makebox(0,0)[c]{$M_{C}$}}
                             \put(43,1){\vector(1,0){3}}
 \put(47,1){\makebox(0,0)[c]{$0$}}

                 \put(21,4.5){\vector(0,-1){2}}
             \put(22,3.5){\makebox(0,0)[c]{$_{f}$}}
             \put(27,4){\makebox(0,0)[c]{$\theta$}}
                 \put(31,4.5){\vector(0,-1){2}}
                 \put(33,3.5){\makebox(0,0)[c]{$_{(\theta,gb)}$}}
                 \put(41,4.5){\vector(0,-1){2}}
                 \put(42,3.5){\makebox(0,0)[c]{$_{g}$}}
                 \put(25,1.5){\makebox(0,0)[c]{$_{i}$}}
                 \put(37,1.5){\makebox(0,0)[c]{$_{\pi}$}}

 \put(14,6){\makebox(0,0)[c]{$0$}}
                              \put(16,6){\vector(1,0){3}}

 \put(21,6){\makebox(0,0)[c]{$A$}}
                              \put(23,6){\vector(1,0){3}}

 \put(31,6){\makebox(0,0)[c]{$B$}}
                              \put(36,6){\vector(1,0){3}}
                              \put(30,5){\vector(-2,-1){7}}
                              \put(37,7){\makebox(0,0)[c]{$_{b}$}}
                              \put(25,7){\makebox(0,0)[c]{$_{a}$}}
 \put(41,6){\makebox(0,0)[c]{$C$}}
                              \put(42,6){\vector(1,0){3}}
 \put(47,6){\makebox(0,0)[c]{$0$}}

\end{picture}

\vskip 10pt

\noindent where $i$ and $\pi$ is canonical injective and canonical projective respectively.
Since the first row in above diagram is $\mathrm{Hom}_R(-,M)$-exact, there exists a morphism $\theta$ such that $f=\theta a$. It further induces a commutative diagram as above.
It follows from Snake Lemma that $(\theta, gb)$ is injective, i.e. $B \in \mathrm{Cogen}M$.

$\Rightarrow$  Since $B \in \mathrm{Cogen}M$ and $M$ is a  costar module,
we have that $0\to B\to M_{0}\to L^{'}\to 0$ with
$M_{0} \in \mathrm{Adp}M$ and $L^{'}\in \mathrm{Cogen}M$. There is a module $M_{0}^{'}$ such that
$M_{0}\oplus M_{0}^{'}=M^{X}$. So we obtain a new short exact sequence
$0\to B\to M^{X}\to L\to0$ with $L=L^{'}\oplus M_{0}^{'}\in \mathrm{Cogen}M$. Consider the pushout of $B\to C$ and
$B\to M^{X}$:

$$  $$

 \setlength{\unitlength}{0.09in}
 \begin{picture}(50,18)

 \put(28,0){\makebox(0,0)[c]{$0$}}
 \put(35,0){\makebox(0,0)[c]{$0$}}

                 \put(28,3.4){\vector(0,-1){2}}
                 \put(31,6){\makebox(0,0)[c]{$1$}}
                 \put(35,3.4){\vector(0,-1){2}}

 \put(28,5){\makebox(0,0)[c]{$L$}}
                             \put(30,5){\vector(1,0){3}}
 \put(35,5){\makebox(0,0)[c]{$L$}}

                 \put(28,9){\vector(0,-1){2}}
                 \put(35,9){\vector(0,-1){2}}

 \put(14,11){\makebox(0,0)[c]{$0$}}
                              \put(16,11){\vector(1,0){3}}
 \put(21,11){\makebox(0,0)[c]{$A$}}
                             \put(23,11){\vector(1,0){3}}
 \put(28,11){\makebox(0,0)[c]{$M^{X}$}}
                             \put(30,11){\vector(1,0){3}}
 \put(35,11){\makebox(0,0)[c]{$N$}}
                             \put(37,11){\vector(1,0){3}}
 \put(41,11){\makebox(0,0)[c]{$0$}}

                 \put(21,14.5){\vector(0,-1){2}}
             \put(22,13.5){\makebox(0,0)[c]{$1$}}
                 \put(28,14.5){\vector(0,-1){2}}
                 \put(35,14.5){\vector(0,-1){2}}

 \put(14,16){\makebox(0,0)[c]{$0$}}
                              \put(16,16){\vector(1,0){3}}

 \put(21,16){\makebox(0,0)[c]{$A$}}
                              \put(23,16){\vector(1,0){3}}

 \put(28,16){\makebox(0,0)[c]{$B$}}
                              \put(30,16){\vector(1,0){3}}
 \put(35,16){\makebox(0,0)[c]{$C$}}
                              \put(37,16){\vector(1,0){3}}
 \put(41,16){\makebox(0,0)[c]{$0$}}

                 \put(28,19.5){\vector(0,-1){2}}
                 \put(35,19.5){\vector(0,-1){2}}

 \put(28,21){\makebox(0,0)[c]{$0$}}

 \put(35,21){\makebox(0,0)[c]{$0$}}

\end{picture}

$$$$

Since the first row and second column are $\mathrm{Hom}_R(-,M)$-exact in above diagram, it is easy
to see that the second row is also $\mathrm{Hom}_R(-,M)$-exact.
Since $A \in \mathrm{Cogen}M$, similar to $B$, there is a short exact sequence
$0\to A\to M^{Y}\to K\to 0$ with $K \in \mathrm{Cogen}M$.
By Lemma \ref{direct}, we have that $M^{Y}\oplus N\cong M^{X}\oplus K$, then it is
easy to see that $N \in \mathrm{Cogen}M$. Consequently,
$C \in \mathrm{Cogen}M$.
\ \hfill $\Box$

\vskip 10pt
Now we give some characterizations of quasi-cotilting modules.

\bg{Pro}\label{LMBAR}%
Let $M$ be an $R$-Module, the following statements are equivalent:

$(1)$ $M$ is a quasi-cotilting module;

$(2)$ $\mathrm{Cogen}M=\mathrm{Copres}M$ and $M$ is Ext-injective in $\mathrm{CogenM}$;

$(3)$ $M$ is a costar module and $\mathrm{Cogen}M$ is a torsion-free class;

$(4)$ $\mathrm{Cogen}M=\mathrm{Fac}(\mathrm{Cogen} M)\bigcap$ $^{\bot_{1}}M$.

\ed{Pro}

\Pf. (1)$\Rightarrow$ (2) By the involved definitions.

(2)$\Rightarrow$ (1),(3) If $M$ is Ext-injective in $\mathrm{Cogen}M$, then it is easy to see that the functor $\mathrm{Hom}_R(-,M)$ preserves the exactness of any short exact sequences in $\mathrm{Cogen}M$. So $M$ is a costar module by the assumption. By Lemma \ref{IARM}, we also have that  $\mathrm{Cogen}M$ is a torsion-free class.

(3)$\Rightarrow$ (4) Clearly, $\mathrm{Cogen}M\subseteq\mathrm{Fac}(\mathrm{Cogen} M)$. To see that $\mathrm{Cogen}M\subseteq {^{\bot_{1}}M}$, we take any extension $0\to M\to N\to L\to 0$ with $L\in\mathrm{Cogen}M$. Since $\mathrm{Cogen}M$ is a torsion-free class, it is closed under extensions. It follows that $N\in\mathrm{Cogen}M$. But $M$ is a costar module, the exact sequence is then $\mathrm{Hom}_R(-,M)$-exact. Thus, we can obtain that the exact sequence is actually split. Consequently we have that $L\in {^{{\bot}_1}M}$ for any $L\in\mathrm{Cogen}M$, i.e.,  $\mathrm{Cogen}M\subseteq {^{\bot_{1}}M}$.

On the other hand, take any $N\in \mathrm{Fac}(\mathrm{Cogen} M)\bigcap$ $^{\bot_{1}}M$.
Then there is a module $L$ in $\mathrm{Cogen} M$ and an epimorphism $f$: $L\to N$. Set $K=\ker f$,
then $K\in \mathrm{Cogen} M$
since $L\in \mathrm{Cogen} M$. Then there is a short exact sequence
$0\to K\to M_{0}\to A\to 0$ with $A \in \mathrm{Cogen}M$ and $M_{0}\in \mathrm{Adp}M$, as $\mathrm{Cogen}M=\mathrm{Copres}M$. Since $N\in$ $^{\bot_{1}}M$,
we have the following communicative diagram:

 \setlength{\unitlength}{0.09in}
 \begin{picture}(50,10)

 \put(14,1){\makebox(0,0)[c]{$0$}}
                              \put(16,1){\vector(1,0){3}}
 \put(21,1){\makebox(0,0)[c]{$K$}}
                             \put(23,1){\vector(1,0){3}}
 \put(28,1){\makebox(0,0)[c]{$M_{0}$}}
                             \put(30,1){\vector(1,0){3}}
 \put(35,1){\makebox(0,0)[c]{$A$}}
                             \put(37,1){\vector(1,0){3}}
 \put(41,1){\makebox(0,0)[c]{$0$}}

                 \put(21,4.5){\vector(0,-1){2}}
             \put(22,3.5){\makebox(0,0)[c]{$1$}}
                 \put(28,4.5){\vector(0,-1){2}}
                 \put(29,3.5){\makebox(0,0)[c]{$_{\alpha}$}}
                 \put(35,4.5){\vector(0,-1){2}}
                 \put(36,3.5){\makebox(0,0)[c]{$_{\beta}$}}

 \put(14,6){\makebox(0,0)[c]{$0$}}
                              \put(16,6){\vector(1,0){3}}

 \put(21,6){\makebox(0,0)[c]{$K$}}
                              \put(23,6){\vector(1,0){3}}

 \put(28,6){\makebox(0,0)[c]{$L$}}
                              \put(30,6){\vector(1,0){3}}
                              \put(32,7){\makebox(0,0)[c]{$_{f}$}}
 \put(35,6){\makebox(0,0)[c]{$N$}}
                              \put(36,6){\vector(1,0){3}}
 \put(41,6){\makebox(0,0)[c]{$0$}}

\end{picture}

By Snake Lemma, we obtain that $\ker\beta=\ker\alpha\in  \mathrm{Cogen} M=\mathrm{Copres}M$
since $L \in \mathrm{Cogen}M$. From the third column in above diagram,
we obtain a new short exact sequence
$0\to \ker\beta\to N\to \mathrm{Im} \beta\to 0$ with $\mathrm{Im} \beta \in \mathrm{Cogen}M$ (since $A \in \mathrm{Cogen}M$).
By the assumption, $\mathrm{Cogen}M$ is a torsion-free class, so  we have
$N \in \mathrm{Cogen}M$. Hence,
$\mathrm{Cogen}M=\mathrm{Fac}(\mathrm{Cogen} M)\bigcap$ $^{\bot_{1}}M$.

(4)$\Rightarrow$ (2) We need only to prove that $\mathrm{Cogen}M\subseteq\mathrm{Copres}M$.
Take any $N\in \mathrm{Cogen}M$ and consider the evaluation map $u$: $N\to M^{X}$ with $X=\mathrm{Hom}_R(N,M)$.
It is easy to verify that $u$ is injective. Thus we have a short exact sequence
$0\to N\to M^{X}\to C\to 0$, so it is enough to prove that $C\in \mathrm{Cogen}M$.
Applying the functor $\mathrm{Hom}_R(-,M)$ to the sequence, we have an exact sequence
$0\to \mathrm{Hom}_R(C,M)\to \mathrm{Hom}_R(M^{X},M)\to^{\gamma} \mathrm{Hom}_R(N,M)\to \mathrm{Ext}^{1}_R(C,M)\to 0$.
It follows from $u$ is the evaluation map that $\gamma$ is surjective. So that $\mathrm{Ext}^{1}_R(C,M)=0$,
i.e. $C\in$ $^{\bot_{1}}M$.
It follows that $C\in\mathrm{Fac}(\mathrm{Cogen} M)\bigcap$ $^{\bot_{1}}M= \mathrm{Cogen}M$ by (4).
\ \hfill $\Box$

\vskip 10pt

The above result suggests the following definition.

\bg{Def}\label{}%
The torsion-free class $\mathcal{T}$ is called a quasi-cotilting class if $\mathcal{T}=\mathrm{Cogen}M$ for some quasi-cotilting module.
\ed{Def}


Let $M$ be an $R$-module. We denoted by $\mathrm{Ann} M$ the ideal of $R$ consisting of all elements $r \in R$ such that $rM=0$.
If $\mathrm{Ann} M=0$, then $M$ is called faithful.

\bg{Lem}\label{faith}%
The following statements are equivalent for an $R$-module $M$.

$(1)$ $M$ is faithful;

$(2)$ $R\in \mathrm{Cogen} M$;

$(3)$ $\mathrm{Proj}R\subseteq \mathrm{Cogen} M$;

$(4)$ $\mathrm{Fac}(\mathrm{Cogen} M)=R{\mathrm{-Mod}}$.

$(5)$ $Q\in \mathrm{Fac}(\mathrm{Adp} M)$.

$(6)$ $Q\in \mathrm{Fac}(\mathrm{Cogen} M)$.

\ed{Lem}

\Pf. (1)$\Leftrightarrow$(2) Note that $\mathrm{Ann} M$ is just the kernel of the evaluation map $R\to M^{\mathrm{Hom}_R(R,M)}$, so the result follows from the universal property of the evaluation map. 

 (2)$\Leftrightarrow$(3) This is followed from the fact that $\mathrm{Cogen}M$ is closed under direct sums and direct summands.
 
  (3)$\Rightarrow$(4) Using the fact that every module is a quotient of a projective module.
  
  (4)$\Rightarrow$(3) If $P$ is any projective module and $\mathrm{Fac}(\mathrm{Cogen} M)=R{\mathrm{-Mod}}$, then there is an exact sequence $0\to P_1\to C\to P\to 0$ with $C\in \mathrm{Cogen}M$. But $P$ is projective implies that the exact sequence is split, it follows that $P$ is a direct summand of $C$ and consequently, $P\in \mathrm{Cogen}M$. Thus, (3) follows. 
  
  (4)$\Rightarrow$(5)$\Rightarrow$(6) Obviously.

  (6)$\Rightarrow$(4)   
Clearly, $\mathrm{Fac}(\mathrm{Cogen} M)\subseteq R$-Mod. On the other hand, since $Q\in \mathrm{Fac}(\mathrm{Cogen} M)$,
there exists an $H \in \mathrm{Cogen} M$ such that $f$: $H\to Q$ is
surjective. For any $L\in R$-Mod,
we have a monomorphism $L\to Q^{X}$ for some $X$, since $Q$ an injective cogenerator.
Then $f^{X}$ is surjective and $Q^{X}\cong H^{X}/K$ with $K=\ker f^{X}$. Consequently,
there exists a submodule $H_{1}$ of $H^{X}$ such that $L\cong H_{1}/K$.
It is easy to see that $H_{1}\in \mathrm{Cogen} M$, so $L\in \mathrm{Fac}(\mathrm{Cogen} M)$
and (4) holds.
  \hfill $\Box$ 

\vskip 10pt

Recall that an $R$-module $M$ is called (1-)cotilting if it satisfies the following three conditions:
(1) the injective dimension of $M$ is not more than $1$, i.e., $\mathrm{id}M \leq 1$;
(2) $\mathrm{Ext}_R^{1}(M^{\lambda},M)=0$ for any set $\lambda$;
(3) There is an exact sequence
$0\to M_{1}\to M_{0}\to Q\to 0$
where $M_{0}, M_{1}\in \mathrm{Adp}M$ and $Q$ is an injective cogenerator. Note that $M$ is cotilting is equivalent  to that $\mathrm{Cogen}M= {^{\bot_{1}}M}$, see for instance \cite{HTTT}. We will freely use these two equivalent definition of cotilting modules.

A torsion-free class $\mathcal{T}$ is called a cotilting class if $\mathcal{T}=\mathrm{Cogen}M$ for some cotilting module $M$.

We have the following characterizations of cotilting modules. Some of them were obtained in \cite{LHX} (in Chinese). For reader's convenience, we include here a complete proof.

\bg{Pro}\label{LQBA}%
Let $M$ be an $R$-module and $Q$ be an injective cogenerator of $R{\mathrm{-Mod}}$, then the following statements are equivalent:

$(1)$ $M$ is a cotilting module;

$(2)$ $M$ is a quasi-cotilting module and $Q\in \mathrm{Fac}(\mathrm{Adp} M)$;

$(3)$ $M$ is a quasi-cotilting module and $\mathrm{Proj}R\subseteq \mathrm{Cogen} M$;

$(4)$ $M$ is a quasi-cotilting module and $\mathrm{Fac}(\mathrm{Cogen} M)=R{\mathrm{-Mod}}$;

$(5)$ $M$ is a  faithful quasi-cotilting module.

$(6)$ $M$ is a quasi-cotilting module and  $\mathrm{Cogen} M$ is a cotilting torsion-free class.

$(7)$ $M$ is a  costar module and $Q\in \mathrm{Fac}(\mathrm{Adp} M)$;

$(8)$ $M$ is a  costar module and $\mathrm{Proj}R\subseteq \mathrm{Cogen} M$;

$(9)$ $M$ is a  costar  module and $\mathrm{Fac}(\mathrm{Cogen} M)=R{\mathrm{-Mod}}$;

$(10)$ $M$ is a  faithful costar module.

$(11)$ $M$ is a costar module and  $\mathrm{Cogen} M$ is a cotilting torsion-free class.

\ed{Pro}

\Pf. 
(2)$\Leftrightarrow$(3)$\Leftrightarrow$(4)$\Leftrightarrow$(5)$\Rightarrow$(10)$\Leftrightarrow$(9)$\Leftrightarrow$(8)$\Leftrightarrow$(7) By Lemma \ref{faith} and the definitions.

(1)$\Rightarrow$(2) Since $M$ is cotilting, we have that and $\mathrm{Cogen}M= {^{\bot_{1}}M}$ and $\mathrm{Fac}(\mathrm{Cogen} M)=R{\mathrm{-Mod}}$ by the above arguments. It follows $\mathrm{Cogen}M=\mathrm{Fac}(\mathrm{Cogen} M)\bigcap {^{\bot_{1}}M}$. Hence, $M$ is quasi-cotilting by Proposition \ref{LMBAR}.

(8)$\Rightarrow$(1) Since $\mathrm{Proj}R \subseteq\mathrm{Cogen} M$, we have that $\mathrm{Fac}(\mathrm{Cogen} M)$=$R$-Mod.
For any $T\in \mathrm{Cogen} M$, there is a short exact sequence $0\to K\to P_{0}\to T\to 0$
with $P_{0}\in \mathrm{Proj}R$. Note that the sequence is indeed in $\mathrm{Cogen} M$, so $\mathrm{Hom}_R(-,M)$ preserves exactness of this short exact sequence, as $M$ is a costar module. It follows that $\mathrm{Ext}^{1}_R(T,M)=0$ since $\mathrm{Ext}^{1}_R(P_0,M)=0$. Thus
$\mathrm{Cogen} M \subseteq$ $^{\bot_{1}}M$ and $\mathrm{Cogen} M$ is a
torsion-free class by Lemma \ref{IARM}. Furthermore, we have that $M$ is a quasi-cotilting module by Proposition \ref{LMBAR}.
Consequently, $\mathrm{Cogen}M=$$^{\bot_{1}}M$ by 
Proposition \ref{LMBAR} again, since $\mathrm{Fac}(\mathrm{Cogen} M)$=$R$-Mod. So $M$ is a cotilting module.

(1)$\Rightarrow$(6)$\Rightarrow$(11) It is obvious now.

(11)$\Rightarrow$(8). If $\mathrm{Cogen} M$ is a cotilting torsion-free class, i.e., $\mathrm{Cogen} M=\mathrm{Cogen} T$ is for some cotilting module $T$, Then $\mathrm{Proj}R\subseteq \mathrm{Cogen} T$ by the above argument. Thus (8) follows.
\ \hfill $\Box$

\mskip\

The above result yields the following characterization of costar modules.

\bg{Pro}\label{TFAE}%
The following statements are equivalent for an $R$-module $M$.

$(1)$ $M$ is a costar $R$-module;

$(2)$ $M$ is a costar $\bar{R}$-module, where $\bar{R}=R/{\mathrm{Ann}}M$;

$(3)$ $M$ is a costar $R/I$-module for any ideal $I$ of $R$ such that $IM=0$;

$(4)$ $M$ is a cotilting $\bar{R}$-module.
\ed{Pro}

\Pf. Note that the category of $R/I$-modules can be identified to the full subcategory 
$\{M \in$ $R$-Mod $|$ $IM=0\}$. Under this identification, it is not difficult to verify that
$\mathrm{Cogen} _{R/I}M=\mathrm{Cogen} _{R}M=\mathrm{Cogen} _{\bar{R}}M$.
Therefore, (1)$\Leftrightarrow$(2)$\Leftrightarrow$(3) is obvious from the definitions.

(2)$\Rightarrow$(4) Note that $M$ is always faithful as an $\bar{R}$-module, so the conclusion follows from
Proposition \ref{LQBA}.

(4)$\Rightarrow$(2) By Proposition \ref{LQBA}.
\hfill $\Box$

\mskip\

A short exact sequence $0\to A\to B\to C\to 0$ is called pure exact if it is ($M\otimes_{R}-$)-exact
for any right $R$-module $M$.
In this case, $C$ is called a pure quotient module of $B$. A module $M$ is called pure injective
if any pure exact sequence is $\mathrm{Hom}_R(-,M)$-exact.

\bg{Lem}\label{TFSA}
The following statements are equivalent:

$(1)$ $M$ is a pure injective $R$-module;

$(2)$ For any $X$, the short exact sequence $0\to M^{(X)}\to M^{X}\to C\to 0$ is $\mathrm{Hom}_R(-,M)$ exact;

$(3)$ $M$ is a pure injective $\bar{R}$-module, where $\bar{R}=R/{\mathrm{Ann}}M$.

\ed{Lem}

\Pf. (1) $\Leftrightarrow$ (2) by \cite[Lemma 2.1]{BCMA}.
Similar to (1) $\Leftrightarrow$ (2), it is easy to see that (2) $\Leftrightarrow$ (3)
since $M^{(X)}$, $M^{X}$ and $C$ are in $\bar{R}$-Mod.
\ \hfill $\Box$

\mskip\

From the above discussion, we can get the following important properties
of quasi-cotilting modules.

\bg{Pro}\label{main}
All costar modules are pure injective and cofinendo. Specially,
all quasi-cotilting modules are pure injective and cofinendo.

\ed{Pro}

\Pf. Let $M$ be a costar $R$-module. We obtain that $M$ is a cotilting $\bar{R}$-module
by Proposition \ref{TFAE}. It follows that $M$ is a pure injective $\bar{R}$-module
since all cotilting modules are pure injective \cite{BCMA}.
Thus $M$ is a pure injective $R$-module by Lemma \ref{TFSA}.

To prove that $M$ is cofinendo, we only need to prove that
$\mathrm{Cogen} M$ is closed under direct sums and pure quotient modules
by Lemma \ref{IACA} and Lemma \ref{TFAEF}.
Obviously, $\mathrm{Cogen} M$ is closed under direct sums.
Suppose the short exact sequence $0\to A\to B\to C\to 0$ is pure exact and $B\in \mathrm{Cogen} M$.
It follows that $A\in \mathrm{Cogen} M$ and that
$0\to A\to B\to C\to 0$ is $\mathrm{Hom}_R(-,M)$ exact from  $M$ is  pure injective. So $C\in \mathrm{Cogen} M$ by
Proposition \ref{LMBAC}, i.e. $\mathrm{Cogen} M$ is closed under pure quotient modules.
\ \hfill $\Box$

\bg{Lem}\label{IACA}\cite[Theorem 2.5]{HJCP}
If a class $\mathcal{A}$ is closed under pure quotient modules, then the following statements are equivalent:

$(1)$ $\mathcal{A}$ is closed under arbitrary direct sums;

$(2)$ $\mathcal{A}$ is a precover class;

$(3)$ $\mathcal{A}$ is a cover class.

\ed{Lem}

\bg{Cor}\label{main Cor}
If $M$ is a costar module, then $\mathrm{Cogen} M$ is a cover class.
In particular, $\mathrm{Cogen} M$ is a cover class for any quasi-cotilting module M.
\ed{Cor}
\Pf. By the proof of Proposition \ref{main}, we obtain that
$\mathrm{Cogen}M$ is closed
under pure quotient and direct sums. Thus $\mathrm{Cogen}M$ is a cover class
by Lemma \ref{IACA}.
\ \hfill $\Box$

\mskip\

%
%

\section{Quasi-cotilting torsion-free class}


\bg{Lem}\label{IMIA}
If $M$ is a quasi-cotilting module, then $\mathrm{Adp} M$=$\ker\mathrm{Ext}^{1}_R(\mathrm{Cogen}M,-)\bigcap \mathrm{Cogen} M$.

\ed{Lem}

\Pf. Suppose that $N \in \mathrm{Adp} M$. Clearly we get $N$
$\in$ $\ker\mathrm{Ext}^{1}_R(\mathrm{Cogen}M,-)\bigcap$ $\mathrm{Cogen} M$, from (3) in Proposition \ref{LMBAR},.
For the inverse inclusion, take any $L\in\ker\mathrm{Ext}^{1}_R(\mathrm{Cogen}M,-)\bigcap \mathrm{Cogen} M$.
Then there is a short exact sequence $0\to L\to M_{0}\to C\to 0$ with $M_{0}\in \mathrm{Adp} M$
and $C\in \mathrm{Cogen} M$ since $L\in\mathrm{Cogen} M=\mathrm{Copres} M$.
Since $\mathrm{Ext}^{1}(C,L)=0$, we can obtain that this exact sequence is split
and hence $L \in \mathrm{Adp} M$.
\ \hfill $\Box$

\bg{Th}\label{char}
Let $Q$ be an injective cogenerator of $R$-Mod. The following statements are equivalent:

$(1)$ $M$ is a quasi-cotilting module;

$(2)$ $M$ is Ext-injective in $\mathrm{Cogen} M$ and there is an exact sequence
$$0\to M_{1}\to M_{0}\to^{\alpha} Q$$
with $M_{0}$ and $M_{1}$ in $\mathrm{Adp} M$ and $\alpha$ a $\mathrm{Cogen} M$-precover.
\ed{Th}

\Pf. (1)$\Rightarrow$(2) By Proposition \ref{main},
$M$ is a cofinendo module. Then there is a morphism $\alpha$: $M_{0}\to Q$ with  $\alpha$ an $\mathrm{Adp} M$-precover
by Lemma \ref{TFAEF}.
It can be shown that  $\alpha$ is also a $\mathrm{Cogen} M$-precover. Indeed,
suppose that $M'\in \mathrm{Cogen} M$, we proof that $f$ can factor through $\alpha$
for any $f$: $M'\to Q$.
There exists a monomorphism $i$: $M' \to M^{X}$ since $M'\in \mathrm{Cogen} M$.
We have a morphism $g$: $M^{X}\to Q$ such that $f=gi$ since $Q$ is injective.
It follows from $\alpha$ is $\mathrm{Adp} M$-precover
that there is a morphism $h$: $M^{X}\to M_{0}$ such that $g=\alpha h$. Thus $f=gi=\alpha hi$.
So $\alpha$ is a $\mathrm{Cogen} M$-precover.
Now set $M_{1}=\ker \alpha$, we only need to prove that $M_{1}\in \mathrm{Adp} M$.
Consider the exact sequence $0\to M_{1}\to M_{0}\to^{\pi} \mathrm{Im} \alpha\to 0$, it is easy to prove that
$\pi$ is also $\mathrm{Cogen} M$-precover by the definition. For any $N\in \mathrm{Cogen} M$,
applying the functor $\mathrm{Hom}_R(N,-)$ to this exact sequence, we have that
$$0\to \mathrm{Hom}_R(N,M_{1})\to \mathrm{Hom}_R(N,M_{0})\to \mathrm{Hom}_R(N,\mathrm{Im} \alpha)\to \mathrm{Ext}^{1}_R(N, M_{1})\to \mathrm{Ext}^{1}_R(N,M_{0})=0.$$
It follows from $\pi$ is $\mathrm{Cogen} M$-precover that $\mathrm{Ext}^{1}_R(N, M_{1})=0$.
Obviously, $M_{1}\in \mathrm{Cogen} M$. Thus $M_{1}\in \mathrm{Adp} M$ by Lemma \ref{IMIA}.

(2)$\Rightarrow$(1) We only need to prove that $\mathrm{Cogen}M=\mathrm{Copres}M$ by Proposition \ref{LMBAR}.
Suppose that $N\in \mathrm{Cogen}M$. Consider the evaluation map $a$: $N\to M^{X}$ with $X=\mathrm{Hom}_R(N,M)$, which is injective since $N\in \mathrm{Cogen}M$. Set $C=\mathrm{coker}$ $a$, next we prove that $C\in \mathrm{Cogen}M$.
There is a monomorphism $f$: $C\to Q^{Y}$ for some $Y$ since $Q$ is an injective cogenerator.
Now consider the following commutative diagram:

 \setlength{\unitlength}{0.09in}
 \begin{picture}(50,10)

 \put(14,1){\makebox(0,0)[c]{$0$}}
                              \put(16,1){\vector(1,0){3}}
 \put(21,1){\makebox(0,0)[c]{$M_{1}^{Y}$}}
                             \put(23,1){\vector(1,0){3}}
                             \put(25,2){\makebox(0,0)[c]{$_{\beta}$}}

 \put(28,1){\makebox(0,0)[c]{$M_{0}^{Y}$}}
                             \put(30,1){\vector(1,0){3}}
                             \put(32,2){\makebox(0,0)[c]{$_{\alpha}$}}

 \put(35,1){\makebox(0,0)[c]{$Q^{Y}$}}

                 \put(21,5.5){\vector(0,-1){3}}
             \put(22,4){\makebox(0,0)[c]{$_{h}$}}
                 \put(28,5.5){\vector(0,-1){3}}
                 \put(29,4){\makebox(0,0)[c]{$_{g}$}}
                 \put(35,5.5){\vector(0,-1){3}}
                 \put(36,4){\makebox(0,0)[c]{$_{f}$}}
                 \put(26.5,4){\makebox(0,0)[c]{$_{s_{1}}$}}
                 \put(33.5,4){\makebox(0,0)[c]{$_{s_{0}}$}}

 \put(14,7){\makebox(0,0)[c]{$0$}}
                              \put(16,7){\vector(1,0){3}}

 \put(21,7){\makebox(0,0)[c]{$N$}}
                              \put(23,7){\vector(1,0){3}}
                               \put(25,8){\makebox(0,0)[c]{$_{a}$}}

 \put(28,7){\makebox(0,0)[c]{$M^{X}$}}
                              \put(30,7){\vector(1,0){3}}
                              \put(32,8){\makebox(0,0)[c]{$_{b}$}}
                              \put(27,6){\vector(-1,-1){4}}
 \put(35,7){\makebox(0,0)[c]{$C$}}
                              \put(36,7){\vector(1,0){3}}
                              \put(34,6){\vector(-1,-1){4}}

 \put(41,7){\makebox(0,0)[c]{$0$}}

\end{picture}

Since $\alpha$ is $\mathrm{Cogen} M$-precover, there is a morphism g such that $\alpha g=fb$, and then we have that
$\beta h=ga$. Since $a$ is evaluation map and $M_{1}^{Y} \in \mathrm{Adp}M$, we have
$s_{1}$ such that $h=s_{1}a$. It is easy to see that $(g-\beta s_{1})a=0$, and then we have
$s_{0}$ such that $g-\beta s_{1}=s_{0}b$.
So $\alpha s_{0}b=\alpha(g-\beta s_{1})=\alpha g=fb$, thus $f= \alpha s_{0}$ since $b$ is surjective.
Since $f$ is injective, $s_{0}$ is injective. Consequently, $C\in \mathrm{Cogen}M$ and $\mathrm{Cogen}M=\mathrm{Copres}M$.
Then the proof is completed.

\ \hfill $\Box$

\mskip\

Let $\mathcal{D}$ be a class of $R$-modules. A module $M \in \mathcal{D}$ is called an injective object of $\mathcal{D}$
if for any short exact sequence  $0\to X\to Y\to Z\to 0$ in $\mathcal{D}$ is $\mathrm{Hom}_R(-,M)$-exact.
A module $M\in \mathcal{D}$ is called a cogenerator of $\mathcal{D}$ if for any $N$ in $\mathcal{D}$, there is a monomorphism
$i$: $N\to M^{X} $ for some set $X$.

\bg{Lem}\label{LTBA}
Let $\mathcal{T}$ be a class of R-modules and $Q$ be an injective cogenerator of $R$-Mod.
Suppose that an exact sequence $0\to A\to B\to^{\alpha} Q\to^{\pi} K\to 0$
satisfies that $\alpha$ is a $\mathcal{T}$-precover and that $A$ is $\mathrm{Ext}$-injective in $\mathcal{T}$.
Denote $\mathcal{D}$=$\{ N\in$ $R$-Mod $\mid \mathrm{Hom}_R(N,\pi)=0\}$.
Then $\mathcal{T}\subseteq \mathcal{D}$ and $M:=\mathrm{Im} \alpha$ is an injective cogenerator in $\mathcal{D}$.
\ed{Lem}

\Pf. For any $C\in \mathcal{T}$ and any morphism $f$: $C\to Q$, there is a morphism $g$: $C\to B$ such that $f=\alpha g$
since $\alpha$ is a $\mathcal{T}$-precover. Thus $\mathrm{Hom}_R(C,\pi)=0$ and $\mathcal{T}\subseteq \mathcal{D}$.

Take any short exact sequence $0\to X\to Y\to Z\to 0$ in $\mathcal{D}$ and any morphism $h$: $X\to M$.
Consider the following commutative diagram:

 \setlength{\unitlength}{0.09in}
 \begin{picture}(50,10)

 \put(14,1){\makebox(0,0)[c]{$0$}}
                              \put(16,1){\vector(1,0){3}}

 \put(21,1){\makebox(0,0)[c]{$M$}}
                             \put(23,1){\vector(1,0){3}}
                             \put(25,2){\makebox(0,0)[c]{$_{b}$}}

 \put(28,1){\makebox(0,0)[c]{$Q$}}
                             \put(30,1){\vector(1,0){3}}

 \put(35,1){\makebox(0,0)[c]{$K$}}
                             \put(37,1){\vector(1,0){3}}

 \put(41,1){\makebox(0,0)[c]{$0$}}

                \put(31.5,2){\makebox(0,0)[c]{$_{\pi}$}}
                 \put(21,4.5){\vector(0,-1){2}}
                \put(22,3.5){\makebox(0,0)[c]{$_{h}$}}
                 \put(28,4.5){\vector(0,-1){2}}
                 \put(29,3.5){\makebox(0,0)[c]{$_{c}$}}

                 \put(25.5,3.5){\makebox(0,0)[c]{$_{\beta}$}}

 \put(14,6){\makebox(0,0)[c]{$0$}}
                              \put(16,6){\vector(1,0){3}}

 \put(21,6){\makebox(0,0)[c]{$X$}}
                              \put(23,6){\vector(1,0){3}}
                              \put(25,7){\makebox(0,0)[c]{$_{a}$}}

 \put(28,6){\makebox(0,0)[c]{$Y$}}
                              \put(30,6){\vector(1,0){3}}
                              \put(26,5){\vector(-1,-1){3.5}}

 \put(35,6){\makebox(0,0)[c]{$Z$}}
                              \put(36,6){\vector(1,0){3}}

 \put(41,6){\makebox(0,0)[c]{$0$}}

\end{picture}

It follows that from $Q$ is injective that there is a morphisms $c$ such that $bh=ca$.
Since $Y\in \mathcal{D}$, we have a morphism $\beta$ such that $c=b\beta$.
Thus $bh=b\beta a$ and $h=\beta a$ since $b$ is injective.
So the exact sequence is $\mathrm{Hom}(-,M)$-exact, i.e. $M$ is injective
in $\mathcal{D}$. For any $N\in \mathcal{D}$, we have a monomorphism $i$: $N\to Q^{I}$ since
$Q$ is an injective cogenerator.
Consider the exact sequence $0\to M^{I}\to Q^{I}\to K^{I}\to 0$. Since $N\in \mathcal{D}$,
we have that $\mathrm{Hom}_R(N,\pi^{I})\cong(\mathrm{Hom}_R(N,\pi))^{I}=0$, i.e., $\pi^{I}i=0$.
So there is a morphism $\gamma$: $N\to M^{I}$ such that $i=b^{I}\gamma$ and $\gamma$ is injective since $i$
is injective. Consequently, $M$ is a cogenerator in $\mathcal{D}$.
\ \hfill $\Box$

\mskip\
The following result is usually called Wakamatsu's lemma, see for instance \cite{Enochs}.

\bg{Lem}\label{Wakamatsu}($\mathrm{Wakamatsu's}$ $\mathrm{lemma}$)
If $\mathcal{T}$ is a class of modules closed under extensions and if
$\varphi$: $T\to M$ is a $\mathcal{T}$-cover, then $\ker\varphi \in T^{\bot_{1}}$.

\ed{Lem}

\bg{Th}\label{quasi}
Let $\mathcal{T}$ be a torsion-free classes in $R$-Mod.
The following statements are equivalent:

$(1)$ $\mathcal{T}$ is quasi-cotilting torsion-free. i.e.
there exists a quasi-cotilting module M such that $\mathcal{T}=\mathrm{Cogen }M$;

$(2)$ $\mathcal{T}$ is a cover class;

$(3)$ For any R-module $N$, there is an exact sequence
$$0\to A\to B\to^{\alpha} N$$
with $\alpha$ a $\mathcal{T}$-precover and A $\mathrm{Ext}$-injective in $\mathcal{T}$.

\ed{Th}

\Pf. (1)$\Rightarrow$ (2) By Corollary \ref{main Cor}.

(2)$\Rightarrow$ (3) For any $R$-module $N$, there is an exact sequence
$$0\to A\to B\to^{\alpha} N$$
with $\alpha$ $\mathcal{T}$-cover by (2). By Wakamatsu's lemma, we have that
$A\in \mathcal{T}^{\bot_{1}}$. Since $B \in \mathcal{T}$ and $\mathcal{T}$ is a torsion-free class,
$A$ is in $\mathcal{T}$. Thus $A$ is $\mathrm{Ext}$-injective in $\mathcal{T}$.

(3)$\Rightarrow$ (1) Take $N=Q$ with $Q$ an injective cogenerator of $R$-Mod, then we have
an exact sequence $0\to A\to B\to^{\alpha} Q$ with
$\alpha$ a $\mathcal{T}$-precover and $A$ $\mathrm{Ext}$-injective in $\mathcal{T}$ by assumption.
Set $M=A\oplus B$. Then $\mathrm{Cogen} M \subseteq \mathcal{T}$ since $\mathcal{T}$ is
a torsion-free class. On the other hand, for any $L\in \mathcal{T}$, we have a monomorphism
$f$: $L \to Q^{X}$. There is a morphism $g$: $L\to B^{X}$ such that $f=\alpha^{X}g$
since $\alpha^{X}$ is clearly a $\mathcal{T}$-precover.  It follows from $f$ is injective that
$g$ is injective. Thus $L\in \mathrm{Cogen} M$ and $\mathrm{Cogen} M=\mathcal{T}$.
It remains to show that $M$ is Ext-injective in $\mathrm{Cogen}M$ by Theorem \ref{char}.
In fact, by assumption, we have to verify this only to $B$.
Take any exact sequence $0\to B\to N\to T\to 0$ with $T\in \mathrm{Cogen} M$.
Then $N$ is in $\mathrm{Cogen} M$ since $\mathrm{Cogen} M=\mathcal{T}$ is a torsion-free class.
Set $C=\mathrm{Im} \alpha$, consider the pushout of $B\to N$ and $B\to C$:

 \setlength{\unitlength}{0.09in}
 \begin{picture}(50,25)

 \put(28,0){\makebox(0,0)[c]{$0$}}
 \put(35,0){\makebox(0,0)[c]{$0$}}

                 \put(28,3.4){\vector(0,-1){2}}

                 \put(35,3.4){\vector(0,-1){2}}

 \put(21,5){\makebox(0,0)[c]{$0$}}
                              \put(23,5){\vector(1,0){3}}
 \put(28,5){\makebox(0,0)[c]{$C$}}
                             \put(30,5){\vector(1,0){3}}
                             \put(31,4){\makebox(0,0)[c]{$_{c}$}}
 \put(35,5){\makebox(0,0)[c]{$W$}}
                             \put(37,5){\vector(1,0){3}}
                             \put(38,4){\makebox(0,0)[c]{$_{e}$}}
 \put(42,5){\makebox(0,0)[c]{$T$}}
                             \put(44,5){\vector(1,0){3}}
 \put(48,5){\makebox(0,0)[c]{$0$}}

                 \put(28,8.5){\vector(0,-1){2}}
                 \put(29,7.5){\makebox(0,0)[c]{$_{\pi}$}}
                 \put(35,8.5){\vector(0,-1){2}}
                 \put(42,8.5){\vector(0,-1){2}}
                 \put(43,8){\makebox(0,0)[c]{$_{1}$}}
                 \put(33,7.5){\makebox(0,0)[c]{$_{s_{1}}$}}
                 \put(36,7.5){\makebox(0,0)[c]{$_{d}$}}
                 \put(40,7.5){\makebox(0,0)[c]{$_{s_{0}}$}}

 \put(21,10){\makebox(0,0)[c]{$0$}}
                              \put(23,10){\vector(1,0){3}}

 \put(28,10){\makebox(0,0)[c]{$B$}}
                              \put(30,10){\vector(1,0){3}}
                               \put(31,11){\makebox(0,0)[c]{$_{a}$}}

 \put(35,10){\makebox(0,0)[c]{$N$}}
                              \put(37,10){\vector(1,0){3}}
                              \put(38,11){\makebox(0,0)[c]{$_{b}$}}
                              \put(38,8.5){\makebox(0,0)[c]{$_{\delta}$}}
                              \put(33,9){\vector(-1,-1){3.5}}
 \put(42,10){\makebox(0,0)[c]{$T$}}
                              \put(44,10){\vector(1,0){3}}
                              \put(40,9.5){\vector(-1,0){3}}
                              \put(40,9){\vector(-1,-1){3.5}}

 \put(48,10){\makebox(0,0)[c]{$0$}}

                 \put(35,13.5){\vector(0,-1){2}}
                 \put(28,13.5){\vector(0,-1){2}}

 \put(28,15){\makebox(0,0)[c]{$A$}}
                              \put(30,15){\vector(1,0){3}}
                               \put(31,16){\makebox(0,0)[c]{$_{1}$}}
 \put(35,15){\makebox(0,0)[c]{$A$}}

                 \put(35,18){\vector(0,-1){2}}
                 \put(28,18){\vector(0,-1){2}}

 \put(28,19){\makebox(0,0)[c]{$0$}}

 \put(35,19){\makebox(0,0)[c]{$0$}}

\end{picture}

By Lemma \ref{LTBA}, the exact sequence $0\to B\to N\to T\to 0$ is $\mathrm{Hom}(-,C)$-exact.
So there is a morphism $s_{1}$ such that $\pi=s_{1}a$.
It is easy to see that $(d-cs_{1})a=0$. So we have $d=cs_{1}+s_{0}b$ in the above diagram for some $s_{0}$.
Since $A$ is Ext-injective in $\mathrm{Cogen}M$, there is a morphism $\delta$: $T\to N$ such that
$s_{0}=d\delta$. So $b=ed=e(s_{0}b+cs_{1})=es_{0}b$ and $es_{0}=1$ since $b$ is surjective.
But $b\delta=ed\delta=es_{0}=1$, thus, $b$ is split.
Consequently, $B$ is Ext-injective in $\mathrm{Cogen}M$.
\ \hfill $\Box$

\mskip\

We say that two quasi-cotilting modules $M_{1}$ and $M_{2}$ are equivalent
if $\mathrm{Adp} M_{1}=\mathrm{Adp} M_{2}$.

\bg{Cor}
There are bijections between

$(1)$ equivalence classes of quasi-cotilting modules;

$(2)$ torsion-free cover classes;

$(3)$ torsion-free classes $\mathcal{T}$ in R-Mod such that every module has a
$\mathcal{T}$-precover with Ext-injective kernel.
\ed{Cor}

\Pf. Let $M_{1}$ and $M_{2}$ be two quasi-cotilting modules, it is easy to prove
that $\mathrm{Adp} M_{1}=\mathrm{Adp} M_{2}$ if and only if
$\mathrm{Cogen} M_{1}=\mathrm{Cogen} M_{2}$ by Lemma \ref{IMIA}.
Now this correspondences can be defined as follows:

(1)$\to$(2): $M\longmapsto \mathrm{Cogen} M$

(2)$\to$(3): $\mathcal{T}\longmapsto \mathcal{T}$

(3)$\to$(1): $\mathcal{T}\longmapsto M$ with $\mathrm{Cogen} M=\mathcal{T}$.
\ \hfill $\Box$
\mskip\

{\small

}

\end{document}